\documentclass[reqno,12pt] {amsart}

\newtheorem{theorem}{Theorem}[section]
\newtheorem{lemma}[theorem]{Lemma}
\newtheorem{corollary}[theorem]{Corollary}

\theoremstyle{definition}

\theoremstyle{remark}
\newtheorem{remark}[theorem]{Remark}

\newcommand{\NSE}{Navier-Stokes equation}

\newcommand{\bR}{\mathbb R}

\newcommand{\mysection}[1]{\section{#1}
\setcounter{equation}{0}}

\begin{document}
\title[Smoothness of solutions to NSEs] {On
the local Smoothness of Solutions of the Navier-Stokes Equations}

\author[H. Dong]{Hongjie Dong}
\address[H. Dong]{127 Vincent Hall,
University of Minnesota, Minneapolis, MN
55455, USA} \email{hjdong@math.umn.edu}

\author[D. Du]{Dapeng Du}
\address[D. Du]{127 Vincent Hall,
University of Minnesota, Minneapolis, MN
55455, USA}
\email{dpdu@math.umn.edu}

\date{2004.8.12}



\begin{abstract}
We consider the Cauchy problem for incompressible Navier-Stokes
equations $u_t+u\nabla_xu-\Delta_xu+\nabla_xp=0,\,\,\text{div}\,
u=0\,\, \text{in}\,\, \bR^d \times \bR^+$ with initial data $a\in
L^d(\bR^d)$, and study in some detail the smoothing effect of the
equation. We prove that for $T<\infty$ and for any positive
integers $n$ and $m$ we have $t^{m+n/2}D^m_tD^{n}_x u\in
L^{d+2}(\bR^d\times (0,T))$, as long as the
$\|u\|_{L^{d+2}_{x,t}(R^d\times(0,T))}$ stays finite.
\end{abstract}

\maketitle

\mysection{Introduction}

We consider the Cauchy problem for incompressible Navier-Stokes
equations in $d$ spatial dimensions:
\begin{align}
                            \label{NSeq}
&u_t+u\nabla_xu-\Delta_xu+\nabla_xp=0,\nonumber \\
&\text{div}\,
u=0\quad\,\, \text{in}\,\, \bR^d \times (0,T),\\
&u(\cdot ,0)  = a \quad\,\,    \text{on}\,\, \bR^d. \nonumber
\end{align}

In a well-known paper \cite{kato} Kato proved that the problem is
locally well-posed for $a\in L^d$. Kato's method is based on
perturbation theory of the Stokes kernel and is different from the
energy methods used in the foundational paper \cite{leray}
by Leray. Kato's results have been generalized by many authors.
(See, for example, \cite{giga}, \cite{iftimie},\cite{tataru},
\cite{Planchon} and \cite{Taylor}). In particular, in
\cite{tataru} Koch and Tataru proved well posedness of (\ref{NSeq})
for $a\in BMO^{-1}$.

By known local regularity theory for Navier-Stokes equations it
can be proved that solutions obtained by Kato's (also known as
mild solutions) immediately become smooth.  However, it seems that
precise smoothness properties of these solutions have not been
studied in detail. In this paper we prove, roughly speaking, that
the smoothing effect of the equations in spaces naturally
appearing in Kato's approach is the same for the heat equation.
Moreover, our method gives a natural and very simple proof of
smoothness of mild solutions. In fact, one can get both existence
of mild solutions and their detailed regularity properties with
essentially the same amount of work that is needed to prove the
existence result. Our main result is that the Cauchy problem
(\ref{NSeq}) for $a\in L^d$ has a unique solution $u$ satisfying
$$
t^{m+n/2}D^m_tD^{n}_x u\in L^{d+2}(\bR^d\times (0,T)),
$$ where
$0<T<\infty$ is a time for which mild solution exists in
$L^{d+2}(\bR^d\times (0,T))$. Very recently, we learnt that Giga
and Sawada \cite{giga2} also studied smoothing properties of mild
solutions. Their approach is similar to ours at a conceptual
level, but at a technical level the results are not the same. In
particular, we obtain extra information about time derivatives.

As most results on mild solutions, our proof is based on a
perturbation argument. The difference is that we work in weighted
function spaces which naturally take into account the smoothing
effect of the Stokes kernel on higher Sobolev norms. Moreover, it
turns out that by suitably organizing the calculations necessary
for obtaining the estimates, the amount of work one has to do is
not much bigger than in proofs which do not give higher
regularity. We remark that if the initial data is small enough,
our proof gives information about the decay of the solution as
$t\to\infty$.

The article is organized as follows. Our main theorem is given in
the following section. An auxiliary result regarding the
regularity of the solution to the heat equation is presented in
section 3. We give the proof of our main theorem in the last
section. And two intermediate lemmas are proved in section 4 and
section 5.

To conclude this Introduction, we explain on some notation used in
what follows: $\bR^d$ is the $d$-dimensional Euclidean space with
a fixed orthonormal basis. A typical point in ${\mathbb R}^d$ is
denoted by $x=(x^1,x^2,...,x^d)$. Various constants are denoted by
$\delta$ and $N$ in general and the expression
$\delta=\delta(\cdots)$ ($N=N(\cdots)$) means that the given
constant $\delta$ ($N$, respectively) depends only on the contents
of the parentheses. $D_t^k$ means taking $k^{\text{th}}$ partial
derivative with respect to $t$ and define
$$
D_x^l:=\frac{\partial^{l_1}}{\partial
x_1^{l_1}}\frac{\partial^{l_2}}{\partial x_2^{l_2}}\cdots
\frac{\partial^{l_d}}{\partial x_d^{l_d}},$$ where multi-index
$l=(l_1,l_2,\cdots,l_d)$. For any vector valued function $v(x,t)$,
$v_j$ means the $j^{\text th}$ component of $v$. We also use the
summation convention over repeated indices. Finally, for any
measurable function $u=u(x,t)$ and any $p,q\in [1,+\infty]$, as
usual we define
$$
\|u(x,t)\|_{L_x^qL_t^p}:=\big\|\|u(x,t)\|_{L_x^q}\big\|_{L_t^p}.
$$

\mysection{The Setting and Main Results}

We are interested in the solution of the incompressible \NSE:
\begin{eqnarray}
u_t+u\nabla u-\Delta_x u+\nabla_x p=0 \label{N1}\\
\text{div}_x\,u=0 \\
u(x,0)=a(x).\label{N3}
\end{eqnarray}

Here, we assume that the initial condition $a(x)$ is in
$L^d(\bR^d)$ and $\text{div}\,a(x)=0$. It is known that
(\ref{N1})-(\ref{N3}) has a unique solution $u(x,t)$ in
$L^{d+2}(\bR^d \times [0,T_u])$, where $T_u$ is the time of
existence of the solution in the class of $L_{x,t}^{d+2}$. Such
$u$ is called the local solution if $T_u<\infty$. Moreover, one
can prove that $u(x,t)$ is smooth in $\bR^d \times (0,T_u]$ by
using a bootstrap argument.

Here comes our main results.

\begin{theorem}
                                                \label{thm1}
Let u(x,t) be the solution of (\ref{N1})-(\ref{N3}) in
$L^{d+2}(\bR^d \times [0,T_u])$ with $T_u<\infty$. Then for any
positive integers $m$ and $n$, we have
$$t^{m+n/2}D_t^m\nabla_x^n u\in L^{d+2}(\bR^d \times
[0,T_u]).$$
\end{theorem}

\begin{remark}
Under more condition that $L^d$ norm of the initial data $a(x)$ is
sufficient small, our method gives some estimate on the decay.
The precise statement is the following: there exists
$\epsilon_0=\epsilon_0(d,m,n)$ such that if
$\|a\|_{L^d(\bR^d)}<\epsilon_0$, then we have
$$t^{m+n/2}D_t^m\nabla_x^nu \in L^{d+2}(\bR^d\times \bR^+)$$.
\end{remark}

\begin{remark}
From \cite{S} we know $u$ is smooth in $\bR^d\times (0,T_u]$. We
also know that $u$ blows up at $t=0$ because the initial data is
only in $L^d$. In some sense, Theorem \ref{thm1} is the
quantitative description of the picture above.
\end{remark}

To prove the main theorem, we need the following two lemmas, which
can be looked as weaker versions of our main theorem.

\begin{lemma}
                                                \label{lemma1}
For any positive integer $n$ and any $p\in [2+d,+\infty]$, $q \in
[d,d+2]$ satisfying the condition $2/p+d/q=1$, we can find
$\delta=\delta(d,n,q)>0$ such that
\begin{equation}
                                                \label{bdd}
\sum_{k=0}^{n}\|t^{k/2}\nabla_x^k u\|_{L_x^q({\mathbb
R}^d)L_t^p(0,\delta)}< +\infty.
\end{equation}
\end{lemma}

\begin{lemma}
                                                \label{lemma2}
For any positive integers $m$ and $n$, and for the same $p,q$ in
Lemma \ref{lemma1}, we can find $\delta=\delta(d,m,n,q)>0$ such
that
\begin{equation}
                                                \label{bdd2}
\sum_{j=0}^m\sum_{k=0}^n\|t^{j+k/2}D_t^j\nabla_x^k u\|
_{L_x^q({\mathbb R}^d)L_t^p(0,\delta)}< +\infty.
\end{equation}
\end{lemma}

\mysection{Auxiliary Result}

As a preliminary, we will first prove a lemma concerning the
regularity property of the solution to the heat equation:
\begin{eqnarray}
u_t(x,t)-\Delta_x u(x,t)=0 \,\,\,&\text{in}& \bR^d \times \bR^+
 \label{heat1}\\
u(x,0)=a(x) \,\,\,&\text{on}& \bR^d. \label{heat2}
\end{eqnarray}

We know that
$U(x,t)=\Gamma(x,t)*a(x)=\int_{\bR^d}{\Gamma(x-y,t)a(y)\,dy}$ is
the solution of (\ref{heat1})-(\ref{heat2}), where
$\Gamma(x,t)=(4\pi t)^{-d/2}e^{-|x|^2/(4t)}$ is the heat kernel.
And, $U$ is smooth in $\bR^d \times \bR^+$.

\begin{lemma}
                                        \label{Uest}
For any positive integer $n$ and any $p\in [2+d,+\infty]$, $q \in
[d,d+2]$ satisfying the condition $2/p+d/q=1$, we have
\begin{equation}
                                                \label{Ubdd}
\|t^{n/2}\nabla_x^n U(x,t)\|_{L_x^q(\bR^d)L_t^p({\mathbb R}^+)}<
N(d,n)\|a(x)\|_{L^d(\bR^d)}.
\end{equation}
\end{lemma}
\begin{proof} This result is known. Here we prefer to
give the proof for the sake of completeness. We adapt idea of
interpolation to prove the lemma with the following three steps.

{\it Step 1.} It is clear that
$$
t^{n/2}\nabla_x^n U(x,t)=[Q_n(x/\sqrt{t})\Gamma(x,t)]*a(x),
$$
where $Q_n$ is a polynomial of degree $n$. Firstly, let's assume
$p=q=2+d$. By Fubini theorem and change of variable, we have
\begin{align}
                                            \label{1ststep}
&\|t^{n/2}\nabla_x^n U(x,t)\|_{L_x^{d+2}(\bR^d)L_t^{d+2}({\mathbb
R}^+)} =\big\|\,\|t^{n/2}\nabla_x^n U(x,t)\|_{L_t^{d+2}
({\mathbb R}^+)}\,\big\|_{L_x^{d+2}(\bR^d)} \nonumber \\
&=\big\|\|\int_{\bR^d}Q_n((x-y)/\sqrt{t})\Gamma(x-y,t)a(y)\,dy\|
_{L_t^{d+2}({\mathbb R}^+)}\big\|_{L_x^{d+2}(\bR^d)} \nonumber
\\
&\leq\big\|\int_{\bR^d}\|Q_n((x-y)/\sqrt{t})\Gamma(x-y,t)\|
_{L_t^{d+2}({\mathbb R}^+)}a(y)\,dy\big\|_{L_x^{d+2}(\bR^d)}
\nonumber \\ &\leq
N(n,d)\big\|\int_{\bR^d}|x-y|^{-d+2/(d+2)}a(y)\,dy\big\|
_{L_x^{d+2}(\bR^d)} \nonumber \\ &\leq
N(d,n)\|a(x)\|_{L^d(\bR^d)}.
\end{align}
The last inequality is due to the fractional integration.

{\it Step 2.} Next, we put $p=+\infty$ and $q=d$. Owing to Young's
inequality, we have
\begin{multline}
                                                    \label{2ndstep}
\|t^{n/2}\nabla_x^nU(x,t)\|_{L_x^d(\bR^d)}=
\|(Q_n(x/\sqrt{t})\Gamma(x,t))*a(x)\|_{L_x^d(\bR^d)}\\
\leq\|(Q_n(x/\sqrt{t})\Gamma(x,t)\|_{L_x^1(\bR^d)}
\|a(x)\|_{L_x^d(\bR^d)}\leq N(n,d)\|a(x)\|_{L_x^d(\bR^d)}.
\end{multline}

{\it Step 3.} For any $d+2<p<+\infty$, due to H\"older's
inequality, we get
\begin{equation}
                                                    \label{holder}
\|\nabla_x^nU(x,t)\|^p_{L_x^q}\leq
\|\nabla_x^nU(x,t)\|^{d+2}_{L_x^{d+2}}
\|\nabla_x^nU(x,t)\|^{p-(d+2)}_{L_x^d}.
\end{equation}

By using (\ref{1ststep})-(\ref{holder}), we have
\begin{align*}
&\int_0^{+\infty}t^{np/2}\|\nabla_x^nU(x,t)\|_{L_x^q(\bR^d)}^p\,dt\\
&\leq
\int_0^{+\infty}t^{np/2}\|\nabla_x^nU(x,t)\|^{d+2}_{L_x^{d+2}}
\|\nabla_x^nU(x,t)\|^{p-(d+2)}_{L_x^d}\,dt\\
&\leq
\big\|t^{n/2}\|\nabla_x^nU(x,t)\|_{L_x^{d}}\big\|_{L^{\infty}_t}^{p-d-2}
\int_0^{+\infty}t^{n(d+2)/2}\|\nabla_x^nU(x,t)\|_{L_x^{d+2}(\bR^d)}^{d+2}\,dt\\
&=\|t^{n/2}\nabla_x^nU(x,t)\|^{p-d-2}_{L_x^dL_t^{\infty}}\,
\|t^{n/2}\nabla_x^nU(x,t)\|^{d+2}_{L_x^{d+2}L_t^{d+2}}\\
&\leq N(d,n)\|a(x)\|_{L_x^d}^p.
\end{align*}
\end{proof}

\begin{remark}
                        \label{remark3.0.1}
Since $U$ satisfies (\ref{heat1})-(\ref{heat2}), for any positive
integers $n$ and $m$, and for the same $p,q$ in Lemma \ref{Uest},
one has
\begin{equation}
                                                \label{Ubdd2}
\|t^{n/2+m}D_t^m \nabla_x^n U(x,t)\|_{L_x^q(\bR^d)L_t^p({\mathbb
R}^+)}< N(d,n+2m)\|a(x)\|_{L^n(\bR^d)}.
\end{equation}
\end{remark}

\mysection{Proof of Lemma \ref{lemma1}}

For any $f\in L^2(\bR^d,\bR^d)$, denote  $f={\mathcal P}f+\nabla
\phi$ to be the Helmholtz decomposition and ${\mathcal P}f\in
L^2(\bR^d,\bR^d)$ is divergence-free. It's known that
\begin{equation}
                                        \label{decomp}
\phi=\partial_{x_j}G*f_j,\quad ({\mathcal
P}f)_i=f_i+\partial_{x_ix_j}G*f_j,
\end{equation}
where $G$ is the Poisson kernel.

Let's look at the following Stokes system:
\begin{eqnarray}
u_t-\Delta_x u+\nabla_x p=f \label{S1}\\
\text{div}_x\,u=0 \\
u(x,0)=a(x).\label{S3}
\end{eqnarray}
The solution of the system (\ref{S1})-(\ref{S3}) can be represented
as the following:
\begin{equation}
                                                    \label{rep}
u(x,t)=U(x,t)+\int_0^t\int_{\bR^d}K_{ij}(y,s)f_j(x-y,t-s)\,dy\,ds,
\end{equation}
where function $U(x,t)$ is defined in the previous section and
$K=(K_{ij}(y,s))_{i,j=1}^d$ is the Stokes kernel. We know that
$$
K_{ij}(y,s)=s^{-d/2}H_{ij}(|y|/\sqrt{s}),
$$ where $H_{ij}$ is a
smooth function on $\bR^d$ and $H_{ij}(y)\sim 1/|y|^d
\,\,\text{as}\,\,y\to \infty.$

Since the solution of (\ref{N1})-(\ref{N3}) is divergence-free, we
have $u_j\partial_{x_j}u_i=\partial_{x_j}(u_ju_i)$. Thus,
(\ref{N1}) is equivalent to the following equation
\begin{equation}
                    \label{N1'}
u_t-\Delta_x u+\nabla_x p=-\partial_j(u_ju).
\end{equation}

Upon using the representation formula (\ref{rep}) with
$-\partial_j(u_ju)$ in place of $f$, we get the following integral
equation:
\begin{equation}
                                    \label{inteq}
u(x,t)=U(x,t)-B(u,u)(x,t),
\end{equation}
where $B(u,v)=\int_0^t\partial_{x_j}K(\cdot,s)\ast
u_jv(\cdot,t-s)\,ds$ is a bilinear form.

Given positive numbers $p$, $q$ satisfying the condition
$2/p+d/q=1$, any $\delta>0$ and nonnegative integers $m$, $n$,
denote
$$
\|u\|_{(p,q,m,n,\delta)}:=\sum_{j=0}^{m}\sum_{k=0}^{n}\|t^{j+k/2}D_t^j\nabla_x^k
u\| _{L_x^q(\bR^d)L_t^p(0,\delta)}.
$$
It is easy to see that for fixed $p$, $q$, $m$, $n$ and $\delta$,
$$X=\{v(x,t) |\, \|v\|_{(p,q,m,n,\delta)}< \infty\}$$
is a Banach space with norm $\|\cdot\|_{(p,q,m,n,\delta)}$. And in
the remaining of this section, we assume $m=0$

\begin{remark}
                                \label{remark4.1}
Because of Lemma \ref{Uest}, $U\in X$ for any $\delta\in
(0,+\infty]$ and any $p\geq 2+d$. By using the property of
Lebesgue integral, for any $\epsilon>0$ there exist $\delta_1>0$
such that for any $\delta\in (0,\delta_1]$, we have
$\|U\|_{(p,q,0,n,\delta)}\leq \epsilon$.
\end{remark}

The following lemma shows that $B: X\times X \to X$ is continuous.

\begin{lemma}
                                        \label{Biscountinous}
Given the assumptions above, for any $u,\,v\in X$, the following
estimate holds.
\begin{equation}
\|B(u,v)\|_{(p,q,0,n,\delta)}\leq
N(q,n,d)\|u\|_{(p,q,0,n,\delta)}\|v\|_{(p,q,0,n,\delta)}.
\end{equation}
\end{lemma}
\begin{proof}
The key idea in this proof is to represent
$t^{k/2}\nabla_x^kB(u,v)$ as a summation of two integrals, both of
which can be estimated. We have
$$
t^{k/2}\nabla_x^kB(u,v)=I_1+I_2,
$$
where
$$
I_1=t^{k/2}\int_0^{t/2}\partial_{x_j}K(\cdot,s)*\nabla_x^k(u_jv(\cdot,t-s))\,ds,
$$
$$
I_2=t^{k/2}\int_{t/2}^t\nabla_x^k
\partial_{x_j}K(\cdot,s)*(u_jv(\cdot,t-s))\,ds.
$$

Now, let's estimate $I_1$ and $I_2$ respectively. Owing to Young's
inequality, for fixed $t>0$, we get
\begin{eqnarray}
                                \label{eq4.9}
\|I_1\|_{L^q_x(\bR^d)}\leq
t^{k/2}\int_0^{t/2}\|\partial_{x_j}K(\cdot,s)*\nabla_x^k(u_jv(\cdot,t-s))\|
_{L^q_x(\bR^d)}\,ds\nonumber\\
\leq
t^{k/2}\int_0^{t/2}\|\partial_{x_j}K(\cdot,s)\|_{L_x^{q/(q-1)}(\bR^d)}
\|\nabla_x^k(u_jv(\cdot,t-s))\| _{L^{q/2}_x(\bR^d)}\,ds.
\end{eqnarray}

Notice that $\partial_{x_j}K(x,s)=s^{-(d+1)/2}H'(|x|/s)$ and
$H'(y)\sim 1/|y|^{d+1}$ as $y\to \infty$. Therefore, we have
\begin{equation}
                                    \label{eq4.10}
\|\partial_{x_j}K(\cdot,s)\|_{L_x^{q/(q-1)}(\bR^d)}\leq N(q,d)
s^{-\frac{q+d}{2q}}.
\end{equation}

Due to Leibniz rule and H\"older's inequality, we obtain
\begin{align}
                                    \label{eq4.11}
\|&\nabla_x^k(u_jv(\cdot,t-s))\| _{L^{q/2}_x(\bR^d)}\\
&=(t-s)^{-k/2}\|(t-s)^{k/2}\nabla_x^k(u_jv(\cdot,t-s))\|
_{L^{q/2}_x(\bR^d)} \nonumber \\
&\leq N(k,d)(t-s)^{-k/2}W(t-s),
\end{align}
where
$$
W(s)=\big(\sum_{l=0}^k\|s^{l/2}\nabla_x^lu(\cdot,s)\|_{L_x^q(\bR^d)}\big)\\
\big(\sum_{l=0}^k\|s^{l/2}\nabla_x^lv(\cdot,s)\|_{L_x^q(\bR^d)}\big)
.$$ Again, by H\"older's inequality,
\begin{equation}
                                    \label{eq4.12}
\|W(s)\|_{L_s^{p/2}(0,\delta)}\leq
N\|u\|_{(p,q,0,n,\delta)}\|v\|_{(p,q,0,n,\delta)}.
\end{equation}
By combining (\ref{eq4.9})-(\ref{eq4.12}) together, we get
\begin{align*}
\|I_1\|_{L^q_x(\bR^d)}&\leq
N(q,k,d)\int_0^{t/2}\big(t/(t-s)\big)^{k/2}s^{-\frac{q+d}{2q}}
W(t-s)\,ds\\
&\leq N(q,k,d)\int_0^{t/2}s^{-\frac{q+d}{2q}} W(t-s)\,ds.
\end{align*}
Note that $-(q+d)/(2q)=-1+1/p$, by using fractional integration, we obtain
\begin{equation}
                                \label{eq4.14}
\|I_1\|_{L_x^q(\bR^d)L_t^p(0,\delta)}\leq
N\|u\|_{(p,q,0,n,\delta)}\|v\|_{(p,q,0,n,\delta)}.
\end{equation}

Similarly, due to Young's inequality, we have
\begin{align}
                                \label{eq4.15}
\|I_2\|_{L^q_x(\bR^d)}&\leq
t^{k/2}\int_{t/2}^t\|\nabla_x^k\partial_{x_j}K(\cdot,s)*(u_jv(\cdot,t-s))\|
_{L^q_x(\bR^d)}\,ds\nonumber \\
&\leq
t^{k/2}\int_{t/2}^t\|\nabla_x^k\partial_{x_j}K(\cdot,s)\|_{L_x^{q/(q-1)}(\bR^d)}
\|(u_jv(\cdot,t-s))\| _{L^{q/2}_x(\bR^d)}\,ds.
\end{align}
H\"older's and the property of Stokes kernel inequality yield
\begin{eqnarray}
\|(u_jv(\cdot,s))\|_{L^{q/2}_x(\bR^d)}\leq
\|u(\cdot,s)\|_{L_x^q(\bR^d)}\|v(\cdot,s)\|_{L_x^q(\bR^d)},\\
\|\nabla_x^k\partial_{x_j}K(\cdot,s)\|_{L_x^{q/(q-1)}(\bR^d)} \leq
N(q,d)s^{-\frac{q(k+1)+d}{2q}}.
\end{eqnarray}
Therefore, we obtain
\begin{align*}
&\|I_2\|_{L^q_x(\bR^d)}\\
&\leq N(q,k,d)\int_{t/2}^t
t^{k/2}s^{-\frac{q(k+1)+d}{2q}}
\|u(\cdot,t-s)\|_{L_x^q(\bR^d)}\|v(\cdot,t-s)\|_{L_x^q(\bR^d)}\,ds\\
&\leq N(q,k,d)\int_0^{t/2}s^{-\frac{q+d}{2q}}
\|u(\cdot,t-s)\|_{L_x^q(\bR^d)}\|v(\cdot,t-s)\|_{L_x^q(\bR^d)}\,ds.
\end{align*}
Owing to fractional integration and H\"older's inequality, we get
\begin{multline}
                                    \label{eq4.20}
\|I_2\|_{L_x^q(\bR^d)L_t^p(0,\delta)}\leq
N(q,k,d)\big\|\|u(\cdot,t)\|_{L_x^q(\bR^d)}\|v(\cdot,t)\|_{L_x^q(\bR^d)}\|
_{L_t^{p/2}(0,\delta)}\\
\leq
N\|u\|_{L_x^q(\bR^d)L_t^p(0,\delta)}\|v\|_{L_x^q(\bR^d)L_t^p(0,\delta)}
\leq N\|u\|_{(p,q,0,n,\delta)}\|v\|_{(p,q,0,n,\delta)}.
\end{multline}

Note that (\ref{eq4.14}) and (\ref{eq4.20}) hold for all
$k=0,1,\cdots,n$. Upon taking summation over $k$, the lemma is
proved.
\end{proof}

Now we are ready to prove Lemma \ref{lemma1}. Let $K=N(q,n,d)$ be
the same as in Lemma \ref{Biscountinous}, which does not depend on
$\delta$. Because of Remark \ref{remark4.1}, we can find and fix
$\delta>0$ such that
\begin{equation}
                                    \label{eq5K}
\|U\|_{(p,q,0,n,\delta)}\leq 1/(5K).
\end{equation}

Denote
$$
{\mathcal S}=\{v\in X\,|\,\|v\|_{(p,q,0,n,\delta)}\leq 1/(3K)\},
$$ which is a closed subset in $X$.
Let's define a map $T:\,X\to X$ as
$$T(v)=U+B(v,v).$$ Due to Lemma
\ref{Biscountinous}, we have $T(\mathcal S)\subset \mathcal S$ and
actually $T$ is a contracting map on $\mathcal S$. Indeed, for any
$v_1,\,v_2\in \mathcal S$,
\begin{align*}
                                        \label{contraction}
&\|Tv_1-Tv_2\|_{(p,q,0,n,\delta)}\\
&=\|B(v_1,v_1)-B(v_2,v_2)\|
_{(p,q,0,n,\delta)}\\
&\leq \|B(v_1,v_1)-B(v_1,v_2)\|_{(p,q,0,n,\delta)}+
\|B(v_1,v_2)-B(v_2,v_2)\|_{(p,q,0,n,\delta)}\\
&= \|B(v_1,v_1-v_2)\|_{(p,q,0,n,\delta)}+
\|B(v_1-v_2,v_2)\|_{(p,q,0,n,\delta)}\\
&\leq K\,\|v_1-v_2\|_{(p,q,0,n,\delta)}(\|v_1\|_{(p,q,0,n,\delta)}+
\|v_2\|_{(p,q,0,n,\delta)})\\
&\leq 2/3\,\|v_1-v_2\|_{(p,q,0,n,\delta)}.
\end{align*}
Owing to the contraction mapping theorem, there exist a unique
$u^*\in \mathcal S$ such that $Tu^*=u^*$. By the uniqueness of the
solution of Navier-Stokes equation, $u^*(x,t)=u(x,t)$ for any
$x\in \bR^d,\,t\in (0,\delta)$, and Lemma \ref{lemma1} is proved.

It is convenient to discuss the following generalized result and
proof of it is provided for the sake of completeness.
\begin{corollary}
                                                \label{cor1}
Lemma \ref{lemma1} remains valid when $n$ is an arbitrary
nonnegative number (not only an integer).
\end{corollary}
\begin{proof}
It is clear that
$$
\sum_{k=0}^{n}\|t^{k/2}\nabla_x^k u\|_{L_x^q({\mathbb
R}^d)L_t^p(0,\delta)}\sim
\big\|t^{n/2}\|u(\cdot,t)\|_{W_x^{n,q}(\bR^d)}\big\|_{L_t^p(0,\delta)}.
$$
Consequently, to prove the corollary, it suffices to prove
\begin{equation}
                                            \label{sobolev}
\big\|t^{m/2}\|u(\cdot,t)\|_{W_x^{m,q}(\bR^d)}\big\|
_{L_t^p(0,\delta)} \leq
\big\|t^{n/2}\|u(\cdot,t)\|_{W_x^{n,q}(\bR^d)}\big\|
_{L_t^p(0,\delta)}^{m/n}
\,\|u\|^{1-m/n}_{L_x^q(\bR^d)L^p_t(0,\delta)},
\end{equation} for any $0<m<n$ and
$\delta<\text{min}\{\delta(n,q),\delta(0,q)\}$.

For any integer $k$, define $P_k$ to be the Littlewood-Paley
projection operator. We have the Littlewood-Paley decomposition
$$
u(\cdot,t)=\sum_{k\in \mathbb{Z}}P_ku(\cdot,t).
$$
And, it's known that for any $f\in W_x^{s,q}$, we have
$$
\|f\|_{W^{s,q}(\bR^d)}\sim \|(\sum_{k\in
\mathbb{Z}}|(1+2^k)^sP_kf|^2)^{1/2}\|_{L^q(\bR^d)}.
$$ with the implicit constant depending on $q$ and $d$.
Thus, by H\"older's inequality, we obtain
\begin{align*}
&\|u(\cdot,t)\|_{W_x^{n,q}(\bR^d)}^{m/n}\,\|u(\cdot,t)\|^{1-m/n}
_{L_x^q(\bR^d)}\\
&\sim \big\|(\sum_{k\in
\mathbb{Z}}|(1+2^k)^nP_ku(\cdot,t)|^2)^{1/2}\big\|
_{L_x^q(\bR^d)}^{m/n}\, \big\|(\sum_{k\in
\mathbb{Z}}|(1+2^k)P_ku(\cdot,t)|^2)^{1/2}\big\|
_{L_x^q(\bR^d)}^{1-m/n}\\
&\geq \big\|\big\{(\sum_{k\in
\mathbb{Z}}|(1+2^k)^nP_ku(\cdot,t)|^2)^{m/n}\,(\sum_{k\in
\mathbb{Z}}|(1+2^k)P_ku(\cdot,t)|^2)^{1-m/n}\big\}^{1/2}\big\|
_{L_x^q(\bR^d)}\\
&\geq \big\|\big\{(\sum_{k\in
\mathbb{Z}}|(1+2^k)^mP_ku(\cdot,t)|^2)\big\}^{1/2}\big\|
_{L_x^q(\bR^d)} \sim \|u(\cdot,t)\|_{W_x^{m,q}(\bR^d)}.
\end{align*}
To obtain (\ref{sobolev}), we only need to apply H\"older's
inequality again.
\end{proof}

\mysection{Proof of Lemma \ref{lemma2}}

We'll use the Sobolev imbedding theorem and H\"older's inequality
to prove the lemma. Firstly, let's prove the following lemma,
which is a further generalization of Lemma \ref{lemma2}.
\begin{lemma}
                                        \label{lemma5.1}
For any positive integers $n_1,n_2,\cdots,n_k$ and any $p\in
[2+d,+\infty]$, $q \in [d,d+2]$ satisfying the condition
$2/p+d/q=1$, we can find $\delta=\delta(d,n_1,n_2,..,n_k,k,q)>0$
such that
\begin{equation}
                                                \label{eqlemma5.1}
\big\|t^{(\sum n_j+k-1)/2}\prod_{j=1}^k\nabla_x^{n_j}
u(x,t)\big\|_{L_x^q({\mathbb R}^d)L_t^p(0,\delta)}< +\infty.
\end{equation}
\end{lemma}
\begin{proof}
Denote $\alpha_j=n_j/2+(k-1)/(2k)$. Corollary \ref{cor1} implies
that
$$
t^{\alpha_j}u(x,t)\in
W_x^{2\alpha_j,\tilde{q}}(\bR^d)L_t^{kp}(0,\delta),
$$
where $\delta=\delta(n_j,k,d,q)$ and $\tilde{q}$ satisfies
$2/(kp)+d/\tilde{q}=1$. By using Sobolev imbedding theorem, we
have
\begin{equation}
\|\nabla_x^{n_j} u(x,t)\|_{L_x^{kq}(\bR^d)} \leq
N\|u(x,t)\|_{W_x^{2\alpha_j,\tilde{q}}(\bR^d)},
\end{equation}
Hence,
$$
t^{\alpha_j}\nabla_x^{n_j}u(x,t)\in
L_x^{kq}(\bR^d)L_t^{kp}(0,\delta)
$$ And we complete the proof by using H\"older's inequality
and putting
$\delta$ the minimum of all such $\delta$'s.
\end{proof}

\begin{remark}
In the sequel, Lemma \ref{lemma5.1} is used only in a special case
when $k=2$.
\end{remark}

To prove Lemma \ref{lemma2}, we will also need the following
lemma.

\begin{lemma}
                                        \label{Piscont}
Recall the definition of $\mathcal P$ in the Helmholtz
decomposition. For any $q>0$, $\mathcal P$ is a continuous linear
operator on $L^q(\bR^d)$. Furthermore, for any differentiable
function $f$ on $\bR^d\times \bR^+$, we have
$$
D_x{\mathcal P}f={\mathcal P}D_xf,\quad D_t{\mathcal P}f={\mathcal
P}D_tf.
$$
\begin{proof}
The second assertion follows immediately from the property of the
convolution. To prove the first part, note that
$$
({\mathcal P}f)_i=f_i+\partial_{x_ix_j}G*f_j,
$$ where $G$ is the Poisson kernel. Therefore,
$$\partial_{x_ix_j}G(x) \sim
1/|x|^n,\quad \partial_{x_k}\partial_{x_ix_j}G(x) \sim
1/|x|^{n+1},$$ and the first part of the lemma is a consequence of
the fractional integration.

Notice that $\mathcal P$ is a Calder\'on-Zygmund operator. Indeed,
the condition of boundedness of $\mathcal P$ on $L^2(\bR^d)$
satisfies automatically because of the definition of Helmholtz
decomposition. So the first assertion of the lemma can also be
interpreted as a conclusion of Calder\'on-Zygmund theorem.
\end{proof}
\end{lemma}

We will prove Lemma \ref{lemma2} by induction on $m$.

In the case of $m=0$, it is just the conclusion of Lemma
\ref{lemma1}. Assume the assertion of Lemma \ref{lemma2} holds for
$m\leq m_0$. Due to the interpolation result as in Corollary
\ref{cor1}, for $m\leq m_0$, Lemma \ref{lemma2} holds for
arbitrary nonnegative number $n$. Let's consider $m=m_0+1$.
Clearly, in order to prove the lemma, it suffices to prove that
for any positive integer $n$, there exists
$\delta=\delta(d,m_0,n,q)>0$ such that
\begin{equation}
                            \label{eq5.3}
\|t^{m_0+1+n/2}D_t^{m_0+1}\nabla_x^n u\| _{L_x^q({\mathbb
R}^d)L_t^p(0,\delta)}< +\infty
\end{equation}

We know that $u$ is smooth in $\bR^d\times (0,T_u)$ and  satisfies
the equation
\begin{equation}
u_t=\Delta_xu-{\mathcal P}(u\nabla_x u).
\end{equation}
Therefore,
$$
t^{m_0+1+n/2}D_t^{m_0+1}\nabla_x^n u=
t^{m_0+1+n/2}D_t^{m_0}\nabla_x^n (\Delta_xu-{\mathcal P}(u\nabla_x
u)).
$$
By inductive assumption, we can find $\delta_1$ such that
\begin{equation}
                \label{eq5.4}
\|t^{m_0+1+n/2}D_t^{m_0}\nabla_x^n \Delta_xu\|_{L_x^q({\mathbb
R}^d)L_t^p(0,\delta_1)}< +\infty
\end{equation}
Upon using similar arguments as in the proof of Lemma
\ref{lemma5.1} with obvious modifications, we can find
$\delta_2=\delta_2(d,m_0,n,q)$ such that for any positive integers
$0\leq j\leq m_0$ and $0\leq k \leq n+1$ we have
\begin{equation}
                                                \label{eq5.5}
\big\|t^{m_0+1+n/2}(D_t^j\nabla_x^{k}
u)(D_t^{m_0-j}\nabla_x^{n+1-k} u)\big\|_{L_x^q({\mathbb
R}^d)L_t^p(0,\delta_2)}< +\infty.
\end{equation}

Finally, by using (\ref{eq5.4})-(\ref{eq5.5}),  Leibniz rule and
Lemma \ref{Piscont}, (\ref{eq5.3}) follows upon setting
$\delta=\text{min}(\delta_1,\delta_2)$. And Lemma \ref{lemma2} is
proved.

\begin{remark}
For small initial data, one has $\delta=\infty$, which immediately
yields the result of Theorem \ref{thm1}. Indeed, as we can see in
the proofs of Lemma \ref{lemma1} and \ref{lemma2}, if
$\|a\|_{L^d}$ is sufficient small (depending on given $m$ and
$n$), (\ref{eq5K}) holds true even for $\delta=\infty$. In some
sense, this implies the decay properties of the solution as
$t\to\infty$.

\end{remark}

\mysection{Proof of Main Theorem}

The following lemma is an immediate consequence of
Calder\'on-Zygmund  theorem.

\begin{lemma}
                \label{lemma6.1}
For any $f\in L^p(\bR^d)$, we can find N=N(d) (does not depend on
s) such that
$$
\|K_{ij}(\cdot,s)*f(\cdot)\|_{L_x^p(\bR^d)}\leq N
\|f(\cdot)\|_{L_x^p(\bR^d)}.
$$
\end{lemma}

It's known that the solution $u$ is smooth in $\bR^d \times
(0,T_u]$. Furthermore, for any positive integers $m$, $n$ and any
positive number $\gamma<T_u$, there exists $N=N(m,n,d,\gamma)$
such that
\begin{equation}
                \label{boundedness}
\|D_t^m\nabla_x^nu(x,t)\|_{L^{\infty}(\bR^d \times
[\gamma,T_u])}\leq N
\end{equation}

By using Lemma \ref{lemma2}, we can find $\delta$ such that
$$
\|t^{m+n/2}D_t^m\nabla_x^n u\|
_{L_x^{d+2}({\mathbb R}^d)L_t^{d+2}(0,\delta)}+
\|t^{m+n/2}D_t^m\nabla_x^n u\|
_{L_x^{d}({\mathbb R}^d)L_t^{\infty}(0,\delta)}< +\infty
$$

Let's put $\gamma=\delta/2$. Denote ${\bar u}(x,t)=u(x,t+\gamma)$
for $t\in [0,T_u-\gamma]$ and $x\in \bR^d$. By the uniqueness of
the solution, ${\bar u}(x,t)$ is the solution of
(\ref{N1})-(\ref{N3}) with initial condition ${\bar
a}(x)=u(x,\gamma)$. Correspondingly, let ${\bar
U}(x,t)=U(x,t+\gamma)$, which is the solution of the heat equation
(\ref{heat1})-(\ref{heat2}) with ${\bar a}(x)$ in place of $a(x)$.

Due to Lemma \ref{Uest}, we have
\begin{eqnarray}
\gamma^{n/2}\|\nabla_x^n{\bar U}(x,t)\|_{L^d_x(\bR^d)L^{\infty}_t
({\mathbb R}^+)}\leq N(d,n)\|a(x)\|_{L^d(\bR^d)},\\
\label{eq6.3} \|\nabla_x^j{\bar
U}(x,t)\|_{W^{n,d}_x(\bR^d)L^{\infty}_t({\mathbb R}^+)}\leq
N(\gamma,d,n)\|a(x)\|_{L^d(\bR^d)}
\end{eqnarray}
Recall that we have the representation formula:
$$
{\bar u}(x,t)={\bar U}(x,t)-\int_0^t\int_{\bR^d}K(y,t-s){\bar u}
\nabla_x {\bar u}(x-y,s)\,dy\,ds.
$$ Upon taking $j$th the spatial derivative, we have
\begin{equation}
\nabla_x^j{\bar u}(x,t)=\nabla_x^j{\bar
U}(x,t)-\int_0^t\int_{\bR^d}K(y,t-s) \nabla_x^j\big({\bar
u}\nabla_x {\bar u}(x-y,s)\big)\,dy\,ds.
\end{equation}

Due to Leibniz rule, $\nabla_x^j\big({\bar u}\nabla_x {\bar u}
(x-y,s)\big)$ is a summation of products $D^{l_1}_x{\bar u}\,
D^{l_2}_x{\bar u}$ , where $l_1+l_2=j+1$ and $0\leq l_1\leq j$.
Now by using triangle inequality, (\ref{eq6.3}), Lemma
\ref{lemma6.1} and boundedness of $D^{l_2}_x{\bar u}$, we get
\begin{align*}
&\|\nabla_x^j{\bar u}(\cdot,t)\|_{L^d_x(\bR^d)}\\
&\leq \|\nabla_x^j{\bar
U}(\cdot,t)\|_{L^d_x(\bR^d)}\\
&\quad +\int_0^t\big\|\int_{\bR^d} K(y,t-s)\nabla_x^j\big({\bar
u}\nabla_x {\bar u}(x-y,s)\big)
\,dy\big\|_{L^d_x(\bR^d)}\,ds\\
&\leq N(\gamma,d,n)\|a(x)\|_{L^d(\bR^d)}\\
&\quad +N(j,d,\gamma)\sum_{k=0}^j
\int_0^t\big\|\int_{\bR^d}K(y,t-s)\nabla_x^k {\bar
u}(x-y,s)\,dy\big\|
_{L^d_x(\bR^d)}\,ds\\
&\leq N(\gamma,d,n)\|a(x)\|_{L^d(\bR^d)}+N(j,d,\gamma)
\sum_{k=0}^j\int_0^t\|\nabla_x^k {\bar
u}(x,s)\|_{L^d_x(\bR^d)}\,ds.
\end{align*}

By taking summation over $j=0,1,\cdots,n$, we obtain
\begin{equation}
                                \label{gronwall}
\|{\bar u}(\cdot,t)\|_{W^{n,d}_x(\bR^d)}\leq
N(\gamma,d,n)\|a(x)\|_{L^d(\bR^d)}+N(\gamma,d,n)\int_0^t\|{\bar
u}(\cdot,s)\|_{W^{n,d}_x(\bR^d)}\,ds.
\end{equation}

Thanks to Gronwall's inequality, (\ref{gronwall}) implies that
$\|{\bar u}(\cdot,t)\|_ {W^{n,d}_x(\bR^d)}$ is bounded for $t\in
[0,T_u-\gamma]$. Consequently,
$\|\nabla_x^nu(\cdot,t)\|_{L^d_x(\bR^d)}$ is bounded for $t\in
[\gamma,T_u]$. Again, by using (\ref{boundedness}), we reach
\begin{equation}
                        \label{eq6.5}
\|\nabla_x^nu(\cdot,\cdot)\|_{L^{d+2}_x(\bR^d)L_t^{\infty}([\gamma,T_u])}
\leq N(d,T_u,n,\|a(x)\|_{L^d(\bR^d)}),
\end{equation} and the theorem is proved in the special case $m=0$.
For the general cases, i.e. $m\neq 0$, the theorem can be proved
on the basis of the special case in the same way as Lemma
\ref{lemma2} is derived from Lemma \ref{lemma1}.

\section*{Acknowledgment}
The authors would like to express their sincere gratitude to Prof.
V. Sverak for pointing out this problem and giving many useful
comments for improvement. D. Du is also grateful to the Institute of
Mathematics of Fudan University for its hospitality.
Part of the work was done there.

\end{document}